\documentclass[12pt]{amsart}
\usepackage{amssymb}
\usepackage{epsf}

\newtheorem{thm}{Theorem}
\newtheorem{prop}[thm]{Proposition}

\newtheorem{definition}[thm]{Definition}

\newtheorem{remark}[thm]{Remark}
\newenvironment{rem}{\begin{remark}\rm}{\end{remark}}
\newtheorem{example}[thm]{Example}
\newenvironment{ex}{\begin{example}\rm}{\end{example}}

\newcommand{\DOT}{{\setlength{\unitlength}{1pt}\begin{picture}(3.5,2)(1,1)
\put(2.3,2){\circle*{2}}\end{picture}}}

\newcommand{\IP}{{\mathcal I}{\mathcal P}}

\newcommand{\QED}{
\setlength{\unitlength}{1.0pt}%
\begin{picture}(7.5,7.5)
\put(5,-5){\rule{2.5pt}{2.5pt}}
\put(2.5,-2.5){\rule{5pt}{2.5pt}}
\put(2.5,0){\rule{2.5pt}{2.5pt}}
\put(0,2.5){\rule{5pt}{2.5pt}}
\end{picture}\vspace{10pt}}

\begin{document}

\title[Hopf algebras and edge-labeled posets]{Hopf algebras and edge-labeled
posets} 

\author{Nantel Bergeron \and Frank Sottile}

\address{Department of Mathematics and Statistics\\
        York University\\
        Toronto, Ontario M3J 1P3\\
	CANADA}
\email[Nantel Bergeron]{bergeron@mathstat.yorku.ca}
\urladdr[Nantel Bergeron]{http://www.math.yorku.ca/bergeron}
\address{Department of Mathematics\\
        University of Wisconsin\\
        Van Vleck Hall\\
        480 Lincoln Drive\\
        Madison, Wisconsin 53706-1388\\
        USA}
\email[Frank Sottile]{sottile@math.wisc.edu}
\urladdr[Frank Sottile]{http://www.math.wisc.edu/\~{}sottile}
\thanks{Journal of Algebra, to appear.  \copyright 1999 Academic Press}
\date{29 April 1999}
\thanks{First author supported in part by NSERC and CRM grants}
\thanks{Second author supported in part by NSERC grant  OGP0170279}
\subjclass{06A07, 16W30, 05E05}
\keywords{edge-labeled poset, Hopf algebra, quasi-symmetric function, 
	incidence algebra}

\begin{abstract}
Given a finite graded poset with labeled Hasse diagram, we construct a
quasi-symmetric generating function for chains whose labels have fixed
descents. 
This is a common generalization of a generating function for
the flag $f$-vector defined by Ehrenborg and of a symmetric function 
associated to certain edge-labeled posets which arose in the theory of
Schubert polynomials.
We show this construction gives a Hopf morphism 
from an incidence Hopf algebra of edge-labeled posets to the Hopf algebra 
of quasi-symmetric functions.
\end{abstract}

\maketitle

\begin{center}
{\em To the memory of Gian-Carlo Rota}
\vspace{15pt}
\end{center}

Joni and Rota~\cite{Joni_Rota}, and later Schmitt~\cite{Schmitt} construct
Hopf algebras from partially ordered sets,
giving a global algebraic framework for studying partially ordered sets.
To every graded partially ordered set, Ehrenborg~\cite{Ehrenborg} defines a
quasi-symmetric generating function 
for its flag $f$-vector.
He shows that this induces a Hopf morphism from the Hopf
algebra of graded posets to the Hopf algebra of quasi-symmetric functions.

Edge-labeled posets are finite graded partially 
ordered sets, the edges of the Hasse diagrams of which are labeled with
integers. 
Following the construction of Stanley's symmetric function~\cite{Stanley84},
we associate with each such  
poset a quasi-symmetric generating function for maximal 
chains whose sequence of edge labels has fixed descents.
We show that this reduces to Ehrenborg's function in an important special
case and induces a Hopf morphism from the Hopf algebra of
edge-labeled posets to the Hopf algebra of quasi-symmetric functions.

While studying structure constants for Schubert polynomials, 
we defined a symmetric function for any edge-labeled poset with a certain 
symmetry~\cite{BS_skew},
giving a unified construction of skew Schur functions, Stanley 
symmetric functions, and skew Schubert functions.
We show thata this symmetric function equals the quasi-symmetric 
generating function defined here.\bigskip

\noindent{\bf 1. Edge-labeled posets. }
A {\em poset} $P$ is a finite partially ordered set with 
maximal element $\hat{1}$ and minimal element $\hat{0}$.
For $x\leq y$ in $P$, let $[x,y]:=\{z\mid x\leq z\leq y\}$.
A poset $P$ is {\em graded} of rank ${\rm rk}(P):=n$ if every maximal chain
has length $n$.
Let $R(P)$ be the set of maximal chains in a poset $P$.
The rank ${\rm rk}(x)$ of $x\in P$ is ${\rm rk}[\hat{0},x]$.

We say that $x\lessdot y$ is a {\em cover} if $[x,y]=\{x,y\}$.
An {\em edge-labeled poset} is a graded poset whose covers
are labeled with integers.
The sequence of labels in a maximal chain is its
{\em word}.
The {\em descent set} $D(\rho)$ of a maximal chain $\rho$ with word
$w_1\cdot w_2\cdots w_n$ ($n= {\rm rk}(P)$) is 
$$
D(\rho)\ =\ \{j\mid w_j>w_{j+1}\}.
$$
For $I,J\subseteq\{1,\ldots,{\rm rk}(P){-}1\}$, define
\begin{eqnarray*}
d_I(P)&=& |\{\rho\in R(P)\mid D(\rho)=I\}|\\
f_J(P)&=& |\{\rho\in R(P)\mid D(\rho)\subseteq J\}|\quad
      =\quad \sum_{I\subseteq J} d_I(P).
\end{eqnarray*}
By inclusion-exclusion, we have
$$
d_I(P)\ =\ \sum_{J\subseteq I} (-1)^{|I-J|} f_J(P).
$$
Ehrenborg and Readdy~\cite{ER95} noted 
$d_I(P)$ is an analog of the rank-selected M\"obius invariant, for 
edge-labeled posets.

We sometimes use compositions $\alpha$ of $n$ in place of subsets $I$ of 
$\{1,\ldots,n{-}1\}$ to index these numbers, and we wish to go back and 
forth between these two indexing schemes.
Given a subset $I=\{i_1<i_2<\cdots<i_k\}$ of $\{1,\ldots,n-1\}$,
define a composition 
$\alpha(I):=(i_1,i_2-i_1,\ldots,n-i_k)$ of $n$.
Likewise, given a composition $\alpha=(\alpha_1,\ldots,\alpha_k)$ of
$n$, define a subset $I(\alpha)$
so that $\alpha(I(\alpha))=\alpha$.
The length, $\ell(\alpha)$,  of $\alpha=(\alpha_1,\ldots,\alpha_k)$
is $k$.
Let $C(n)$ be the set of compositions of $n$.
\bigskip

\noindent{\bf 2. Quasi-symmetric
functions. }
Gelfand {\em et al.}~\cite{GKal} define the graded Hopf algebra 
$NC_\DOT$ of non-commutative symmetric functions to be the free associative
algebra with one generator $S_i$ of degree $i$ for each $i=1,2,\ldots$. 
The graded Hopf dual of $NC_\DOT$ is the algebra ${\mathcal Q}_\DOT$
of quasi-symmetric functions~\cite{Ges}, which consists of all formal
power series of bounded degree in commuting indeterminates $x_1,x_2,\ldots$
which are quasi-symmetric:
the coefficient of 
$x_{i_1}^{\alpha_1}x_{i_2}^{\alpha_2}\cdots x_{i_k}^{\alpha_k}$
depends only upon $\alpha$ and not on $i_1,\ldots,i_k$, if
$i_1<i_2<\cdots<i_k$.
Thus ${\mathcal Q}_\DOT$ has a basis of {\em monomial 
quasi-symmetric functions} $M_\alpha$ defined by 
$$
M_\alpha\ :=\ \sum_{i_1<i_2<\cdots<i_k} 
x_{i_1}^{\alpha_1}x_{i_2}^{\alpha_2}\cdots x_{i_k}^{\alpha_k}.
$$

The basis of $NC_\DOT$ dual to the $M_\alpha$ are the {\it quasi-Schur
functions} $S^{\alpha}:= S_{\alpha_1}S_{\alpha_2}\cdots S_{\alpha_k}$, where
$\ell(\alpha)=k$.
Thus, for any edge-labeled poset $P$, the linear map
$\psi_P: NC_\DOT \rightarrow {\mathbb Z}$ given by
$$
\psi_P(S^\alpha)\ =\ \left\{ 
\begin{array}{lcl}
f_{I(\alpha)}(P)&&\mbox{if }\alpha \in C({\rm rk}(P))\\
0   &\quad&\mbox{otherwise}
\end{array}\right.
$$
defines a quasi-symmetric function $F_P$. 
It follows that if $P$ is an edge-labeled poset of rank $n$, then 
\begin{equation}\label{eq:third_def}
F_P\ =\ \sum_{\alpha\in C({\rm rk}(P))} f_{I(\alpha)}(P) M_{\alpha}.
\end{equation}
This function is our main object of study.
\medskip

Another basis of ${\mathcal Q}_\DOT$ is the  {\em fundamental
quasi-symmetric functions} $F_{I,n}$. 
For any subset $I$ of $\{1,\ldots,n{-}1\}$, define 
$$
F_{I,n}\ :=\ 
\sum_{\stackrel{\mbox{\scriptsize $j_1\leq j_2\leq\cdots\leq j_n$}}%
{i\in I \Rightarrow j_i<j_{i+1}}}
x_{j_1} x_{j_2}\cdots x_{j_n}.
$$

One checks that
\begin{equation}\label{eq:reciprocity}
\begin{array}{rcl}
F_{I,n}&=&{\displaystyle \sum_{I\subseteq J\subseteq\{1,\ldots,n-1\}}
M_{\alpha(J)}}\\ 
M_\alpha &=&{\displaystyle \sum_{I(\alpha)\subseteq J}
(-1)^{|J-I(\alpha)|}F_{J,n}} \rule{0pt}{18pt}
\end{array}
\end{equation}

The {\it ribbon Schur functions} $R_\alpha$ form a basis of $NC_\DOT$
dual to the $F_{I,n}$, and there is a change of basis between the $S^\alpha$
and the $R_\alpha$ analogous to (\ref{eq:reciprocity}).
The expressions relating $f_I(P)$ to $d_I(P)$ and 
$F_{I,n}$ to $M_\alpha$ give
\begin{eqnarray}
F_P&=& \sum_{I\subseteq\{1,\ldots,n{-}1\}} d_I(P) F_{I,n}
                      \label{eq:first_def}\\ 
&=& \sum_{\rho\in R(P)} F_{D(\rho),n}. \label{eq:second_def}
\end{eqnarray}
This last expression shows that $F_P$ is a generalization of Stanley's
(quasi-)symmetric function $F_w$~\cite[Equation (1)]{Stanley84},
introduced to study reduced decompositions of elements $w$ of the symmetric 
group.
To see this, let $P$ be the interval $[1,w]$ in the weak order on the 
symmetric group, with the label of a cover $u\lessdot v$ the integer $i$,
where $(i,i+1)=v u^{-1}$.
Then $R([1,w])$ is the set of reduced decompositions of $w$,
and our definition~(\ref{eq:second_def}) for $F_{[1,w]}$ coincides with
Stanley's definition of $F_w$.
\bigskip

\noindent{\bf 3. Incidence Hopf algebras. }
See~\cite{Montgomery,Sweedler} for more on Hopf
algebras.
Let ${\mathcal P}$ be a class of graded posets 
closed under taking subintervals and products.
The (reduced) {\em incidence coalgebra}~\cite{Joni_Rota,Schmitt} $\IP$ of
${\mathcal P}$ 
is the graded free abelian group generated by isomorphism classes
of posets in ${\mathcal P}$ with grading induced by the rank of a poset
and coproduct by 
$$
\Delta(P)\ =\ 
\sum_{x\in P}\; [\hat{0},x]\otimes [x,\hat{1}].
$$
The augmentation is given by projecting onto the degree 0 
component.
The product of posets induces an algebra structure on $\IP$ with identity
the class of a one element poset.

We say that edge-labeled posets $P$ and $Q$ are {\em label-equivalent}
if  there is an isomorphism $P\stackrel{\sim}{\longrightarrow}Q$
preserving the numbers $f_I$  of subintervals.
A map preserving the relative order of the edge labels is such a function,
but there are others.
Suppose now that  ${\mathcal P}$ is a class of edge-labeled posets.
The {\em incidence coalgebra} $\IP$ of
${\mathcal P}$ is the graded free abelian group on 
label-equivalence classes in  ${\mathcal P}$,
with coproduct and augmentation as before.

To define an algebra structure on $\IP$, we first
form the product $P\times Q$ of edge-labeled posets $P$ and $Q$.
Recall that a cover $(p,q)\lessdot (p',q')$ in $P\times Q$ has one of two
forms: either $p=p'$ and $q\lessdot q'$ is a cover in $Q$, or else 
$p\lessdot p'$ is a cover in $P$ and $q=q'$.
Label a cover $(p,q)\lessdot(p',q')$ in $P\times Q$ by the label of the
corresponding cover in $P$ or $Q$.

\begin{prop}[Lemma~3.9 of~\cite{BS_skew}]\label{prop:el-product}
Suppose that $P$ and $Q$ are edge-labeled posets with distinct sets of
edge labels.
Then for any composition $\alpha$ of ${\rm rk}(P)+{\rm rk}(Q)$, 
$$
f_{I(\alpha)}(P\times Q)\ =\ 
\sum_{\beta+\gamma=\alpha} f_{I(\beta)}(P) \cdot f_{I(\gamma)}(Q),
$$
where $\beta$ ranges over compositions of ${\rm rk}(P)$ and 
$\gamma$ over compositions of ${\rm rk}(Q)$, and addition of compositions is
component-wise.
\end{prop}

Let $x,y\in \IP$ be label-equivalence classes of
edge-labeled posets. 
Then $xy$ is the equivalence class with representative $P\times Q$,
where $P$ is a representative of $x$, $Q$ is a representative of $y$, and
$P,Q$ have disjoint sets of edge labels.
This product is independent of choices, by 
Proposition~\ref{prop:el-product}. 
It is also commutative and compatible with
the coproduct, so $\IP$ is a graded bialgebra
and hence has a unique antipode~\cite[Lemma 2.1]{Ehrenborg}.
We summarize these facts.

\begin{thm}
Let ${\mathcal P}$ be a class of edge-labeled posets closed under 
taking subintervals and products.
Then, with the above definitions, $\IP$ is a commutative graded 
Hopf algebra.
\end{thm}

We give our main theorem.

\begin{thm}\label{thm:main}
Let  ${\mathcal P}$ be a class of edge-labeled posets closed under 
subintervals and products.
Then the map $\Phi:\IP\rightarrow {\mathcal Q}_\DOT$ induced by
$$
P\in{\mathcal P}\ \longmapsto\ F_P \in {\mathcal Q}_\DOT
$$
is a morphism of graded Hopf algebras.
\end{thm}

\noindent{\bf Proof. }
The expression~(\ref{eq:second_def}) shows that $F_P$ is a generalization
of Stanley's symmetric function $F_w$.
In fact, the proof~\cite[Theorem 3.4]{Stanley84} that 
$F_{w\times v} = F_w \cdot F_v$ also shows the corresponding fact for $F_P$:
If $P,Q$ are edge-labeled posets with disjoint sets of edge labels, then 
$F_{P\times Q}=F_P\cdot F_Q$.
Thus $\Phi$ is an algebra morphism.

We show it is a coalgebra morphism.
For a composition $\alpha=(\alpha_1,\ldots,\alpha_k)$ and 
integer $0\leq j\leq k (=\ell(\alpha))$, 
define (possibly empty) compositions $\alpha_{\leq j}$ and $\alpha_{>j}$:
\begin{eqnarray*}
\alpha_{\leq j}&:=& (\alpha_1,\ldots,\alpha_j)\\
\alpha_{>j}&:=& (\alpha_{j+1},\ldots,\alpha_k)
\end{eqnarray*}
The coalgebra structure on ${\mathcal Q}_\DOT$ is given by
$$
\Delta M_\alpha\ =\ \sum_{j=0}^{\ell(\alpha)} 
M_{\alpha_{\leq j}}\otimes M_{\alpha_{>j}}.
$$
For an edge-labeled poset $P$ and composition $\alpha$ of ${\rm rk}(P)$, let
$f_{\alpha}(P)=f_{I(\alpha)}(P)$.
Then, for any $1\leq j\leq k$, the following identity is straightforward.
\begin{equation}\label{eq:d-identity}
f_\alpha(P)\ =\ 
\sum_{\stackrel{\mbox{\scriptsize$x\in P$}}{{\rm rk}(x)=I(\alpha)_j}}
f_{\alpha_{\leq j}}[\hat{0},x] \cdot f_{\alpha_{>j}}[x,\hat{1}].
\end{equation}

Using equations~(\ref{eq:third_def}) and~(\ref{eq:d-identity}), we have
\begin{eqnarray*}
\Delta F_P
&=& \sum_{\alpha\in C({\rm rk}(P))} f_\alpha(P) \Delta M_\alpha\quad\ 
=\quad \sum_{\alpha\in C({\rm rk}(P))} f_\alpha(P) \sum_{j=0}^{\ell(\alpha)}
                M_{\alpha_{\leq j}}\otimes M_{\alpha_{>j}}\\
&=& \sum_{\alpha\in C({\rm rk}(P))}\sum_{j=0}^{\ell(\alpha)}
    \sum_{\stackrel{\mbox{\scriptsize$x\in P$}}{{\rm rk}(x)=I(\alpha)_j}}
    f_{\alpha_{\leq j}}[\hat{0},x] M_{\alpha_{\leq j}}\otimes 
    f_{\alpha_{>j}}[x,\hat{1}] M_{\alpha_{>j}}\\
&=& \sum_{x\in P}\left(\sum_{\beta\in C({\rm rk}([\hat{0},x]))} 
f_\beta[\hat{0},x] M_\beta\right)
    \otimes\left(\sum_{\gamma\in C({\rm rk}([x,\hat{1}]))} 
    f_\gamma [x,\hat{1}] M_\gamma\right),
\end{eqnarray*}
which we recognize as $F_{\Delta P}$.
\QED\medskip

\begin{ex}
A Boolean poset is the poset of subsets of a finite set of integers
in which a cover $X\lessdot Y$ is labeled by the integer $X\setminus Y$.
The one element chain $x:= (\hat{0}<\hat{1})$ is the unique primitive
element in any non-trivial (reduced) incidence Hopf algebra of edge-labeled
posets. 
(All labelings of $x$ are equivalent.)
This primitive element generates the commutative subalgebra 
${\mathbb Z}[x]$, which is the incidence Hopf algebra for the class 
${\mathcal B}$ of Boolean posets (algebras) with a standard labeling for
a lattice of order ideals, as the Boolean poset of subsets of
$\{1,2,\ldots,n\}$ is the lattice of order ideals of the antichain
$\{1,2,\ldots,n\}$~\cite[Example 3.13.3]{Stanley_EC1}. 
Moreover, the map 
$\Phi: {\mathcal I}{\mathcal B} (= {\mathbb Z}[x])
\rightarrow {\mathcal Q}_\DOT$ 
is an isomorphism onto the subalgebra
generated by $h_1=F_x$, which is a subalgebra of symmetric functions.
\end{ex}

\noindent{\bf 4. Rank-selected posets. }
Let $P$ be a graded poset and $I$ be a subset of 
$\{1,\ldots,{\rm rk}(P){-}1\}$.
The {\em rank-selected poset} $P(I)$ is the induced subposet of $P$
consisting of all elements of $P$ with rank in $I$, together with $\hat{0}$
and $\hat{1}$. 
Set $\varphi_I(P)$ to be the number of maximal chains in $P(I)$.
These numbers $\varphi_I(P)$ constitute the {\em flag $f$-vector} of $P$.
Ehrenborg's quasi-symmetric generating function $E_P$ for the flag $f$-vector 
satisfies
$$
E_P\ =\   \sum_I \varphi_I(P) M_{\alpha(I)}.
$$

An edge-labeled poset $P$ is {\em $R$-labeled} if every
interval  has a unique increasing chain.
For these posets, $\varphi_I(P)=f_I(P)$, and so $E_P=F_P$.
Similarly, the numbers $d_I(P)$ are the rank-selected M\"obius invariant
for $R$-labeled posets $P$~\cite[Section 3.13]{Stanley_EC1}, and the 
$\eta$ and $\nu$ functions of Ehrenborg-Readdy~\cite{ER95} (for 
edge-labeled posets) reduce to the zeta and M\"obius functions 
for $R$-labeled posets.

More generally, we regard $f_I(P)$ as an extension of the notion of 
flag $f$-vector. 
Suppose $P$ is an edge-labeled poset and 
$I\subseteq\{1,\ldots,{\rm rk}(P){-}1\}$.
Let $P(I)_{\rm wt}$ be the rank selected poset as before, but with 
every cover $x\lessdot y$ in $P(I)_{\rm wt}$ weighted by the number of 
chains with increasing labels in the interval $[x,y]$ of $P$.
A maximal chain in $P(I)_{\rm wt}$ has weight given by the product of the
weights of its covers.
Then $f_I(P)$ counts these weighted maximal chains of $P(I)_{\rm wt}$, 
and so $F_P$ is a weighted version of $E_P$.

Another connection between these theories is given by
Stanley~\cite{Stanley_EJC}
and concerns a relative version of $E_P$ and the flag $f$-vector.
Let $\Gamma$ be a set of (not necessarily maximal) chains of a
poset $P$ that is closed under taking subchains.
The relative flag $f$-vector $\varphi_I(P/\Gamma)$ 
counts chains in the rank selected poset $P(I)$ that are not in $\Gamma$,
and $E_{P/\Gamma}$ is the quasi-symmetric generating function for
$\varphi_I(P/\Gamma)$. 

An edge-labeled poset is {\it relative $R$-labeled\,} if each interval has at
most one increasing chain, and all subintervals of an interval with an
increasing chain also have an increasing chain.
If $P$ is relative $R$-labeled, and $\Gamma$ is the set of chains
$\hat{0}=t_0<t_1<\cdots<t_r=\hat{1}$ for which there is an $i$ where the 
interval $[t_{i-1},t_i]$ does not have an increasing chain, then Stanley
shows that $\varphi_I(P/\Gamma)=f_I(P)$, so that $E_{P/\Gamma}=F_P$.

Interestingly, the labeled posets whose study led us to consider $F_P$
are all relative $R$-labeled.
These are  intervals in the weak order~\cite{Stanley84}, the $k$-Bruhat
order~\cite{BS98}, and the Grassmannian Bruhat order~\cite{BS_monoid}, all
on the symmetric group.\bigskip

\noindent{\bf 5. Symmetric edge-labeled posets. }
For more on symmetric functions, see~\cite{Macdonald95}.
For a composition $\alpha$, let $\lambda(\alpha)$ be the partition obtained
by listing the components of $\alpha$ in decreasing order.
For a partition $\mu$, the {\em monomial symmetric function} $m_\mu$ is
$$
m_\mu \ :=\ \sum_{\alpha\ :\ \lambda(\alpha)=\mu} M_\alpha.
$$
These form a basis for the algebra of symmetric functions.
From Equation~(\ref{eq:third_def}), we deduce the following fact.

\begin{thm}
The function $F_P$ is symmetric if and only if for every 
composition $\alpha$ of\/ ${\rm rk}(P)$, the number $f_\alpha(P)$ 
depends only upon $\lambda(\alpha)$.
\end{thm}

An edge-labeled poset is {\em symmetric} if  
$f_\alpha(P)$ depends only upon $\lambda(\alpha)$.
Symmetric posets arose in the study of Schubert
polynomials~\cite{BS_skew}, where we defined a symmetric function $S_P$ for
each symmetric poset.
This provided a common definition of Stanley symmetric functions, 
skew Schur functions, and skew Schubert functions.
For these, the posets were intervals in, respectively, the weak order on
the symmetric group, Young's lattice, and the Grassmannian Bruhat
order~\cite{BS_monoid}.
The labeling for the weak order was described in Section 2.
In Young's lattice, a cover $\mu\lessdot\lambda$ has a unique index $i$ with
$\mu_i<\lambda_i$, and we label that cover with the the integer
$i-\lambda_i$.
The Grassmannian Bruhat order is a common generalization of both of these
labeled posets, and we refer the reader to~\cite{BS98} for details.

We will show $S_P$ is just the function $F_P$.
A quasi-symmetric generating function construction of skew Schur 
functions was given in~\cite{Ges}, which is essentially 
the same as given here.
While Gessel uses a poset labeling different from that used in~\cite{BS_skew}, 
the two are label-equivalent in a strong sense---a maximal chain in either
labeling has the same descent set.\medskip

The algebra $\Lambda_\DOT$ of symmetric functions has several distinguished
bases besides the $m_\lambda$.
These include the complete symmetric functions $h_\lambda$ and the Schur
functions $S_\lambda$.
These bases are related by the Cauchy formula, an element in the 
graded completion of $\bigoplus_n\Lambda_n(x)\otimes\Lambda_n(y)$:
$$
\prod_{i,j}(1-x_iy_j)^{-1}\ =\ \sum_\lambda h_\lambda(x)m_\lambda(y)
\ =\ \sum_\lambda S_\lambda(x)S_\lambda(y)
$$

Suppose $P$ is a symmetric edge-labeled poset.
Define a linear map $\chi_P:\Lambda_\DOT\rightarrow {\mathbb Z}$ by
$$
\chi_P(h_\mu)\ =\ \left\{\begin{array}{lcl}
f_\mu(P)&&\mbox{if $\mu$ is a partition of }{\rm rk}(P),\\
0&\quad&\mbox{otherwise}.\end{array}\right.
$$
We define $c^P_\lambda := \chi_P(S_\lambda)$, which 
generalizes the Littlewood-Richardson numbers, as
$c^P_\lambda=c^\nu_{\mu,\lambda}$ when $P$ is the interval $[\mu,\nu]$ in
Young's lattice with a natural labeling of covers~\cite{BS_skew}.

Since $\Lambda_\DOT$ is a self-dual Hopf algebra (with $\{h_\mu,m_\mu\}$ and
$\{S_\mu,S_\mu\}$ pairs of dual bases), $\chi_P$ gives a symmetric
function $S_P$.
From the Cauchy formula, we see that 
\begin{eqnarray*}
S_P&=& \chi_P\otimes 1_{\Lambda_\DOT(y)}
\left(\prod_{i,j}(1-x_iy_j)^{-1}\right)\\
&=& \sum_{\lambda\vdash {\rm rk}(P)} f_\lambda(P)\, m_\lambda\\
&=& \sum_{\lambda\vdash {\rm rk}(P)} c^P_\lambda  S_\lambda.
\end{eqnarray*}

\begin{thm}
Let $P$ be a symmetric edge-labeled poset.
Then $S_P=F_P$.
\end{thm}
\noindent{\bf Proof. }
\begin{eqnarray*}
F_P&=&\sum_{\alpha\in C({\rm rk}(P))} f_\alpha(P) M_\alpha\\
&=& \sum_{\mu\vdash {\rm rk}(P)} f_\mu(P) 
      \sum_{\alpha\ :\ \lambda(\alpha)=\mu}M_\alpha\\
&=& \sum_{\mu\vdash {\rm rk}(P)} 
   f_\mu(P) m_\mu \quad = \quad S_P.\qquad \QED
\end{eqnarray*}

\begin{rem}
The definition of $F_P$ in terms of the linear map $\psi_P$ 
(Section 2)  mimics the Cauchy identity construction of $S_P$ above.
The Cauchy identity for $NC_\DOT$ and ${\mathcal Q}_\DOT$ is an element
in the graded completion of 
$\bigoplus_n NC_n\otimes {\mathcal Q}_n$~\cite[Section 6]{GKal}: 
$$
\sum_\alpha R_\alpha\otimes F_{I(\alpha)}\ =\ 
\sum_\alpha S^{\alpha}\otimes M_\alpha.
$$
Thus $F_P$ is just $\psi_P\otimes 1_{{\mathcal Q}_\DOT}$ 
applied to this element.
Here $S^{\alpha}$ is the analog of the homogeneous symmetric function
and $\psi_P(S^{\alpha})=f_{I(\alpha)}(P)$.
\end{rem}

\begin{rem}
In many cases when $F_P$ is symmetric, the symmetric function $F_P$ is known
to be the Frobenius characteristic of a representation of the symmetric
group ${\mathcal S}_{{\rm rk}(P)}$ on the linear span of maximal chains of
$P$.
For example, if $P$ is the Boolean poset of subsets of $[n]$, then
$F_P=(h_1)^n$, which is the Frobenius characteristic of the right regular 
representation of ${\mathcal S}_n$.
This is the action of ${\mathcal S}_n$ on maximal chains induced by
permuting the factors of $P=(\hat{0}<\hat{1})^n$.

Similarly, $F_P$ is the Frobenius characteristic of a representation 
if $P$ is an interval in either the weak order on the symmetric
group~\cite{Kraskiewicz} or Young's lattice~\cite{Sagan}.
If $P$ is an interval in either the $k$-Bruhat order or the Grassmannian
Bruhat order, then $F_P$ is known to be Schur-positive, by geometry.
For such intervals $P$, it would be interesting to construct a 
${\mathcal S}_{{\rm rk}(P)}$-representation on the linear span of the maximal
chains of $P$ with Frobenius characteristic $F_P$.
Considering rank 3 intervals in these orders shows that this representation
cannot arise from a permutation action of ${\mathcal S}_{{\rm rk}(P)}$ on
the maximal chains of $P$.
\end{rem}

\section*{Acknowledgments}
We are indebted to Richard Stanley, who suggested studying chains with fixed
descents using quasi-symmetric functions, and pointed out the connection to
relative $R$-labelings.
We also thank Richard Ehrenborg for helpful conversations and 
Sarah Witherspoon
whose Hopf-algebraic advice was crucial at the start of this project.

\end{document}